\documentclass[10pt,a4paper]{article}
\usepackage[utf8]{inputenc}
\usepackage[english]{babel}
\usepackage{amsmath}
\usepackage{amsfonts}
\usepackage{amssymb}
\usepackage{amsthm}
\usepackage{graphicx}
\usepackage{cite}
\usepackage{algorithmic}
\usepackage{textcomp}
\usepackage{xcolor}
\usepackage{paralist}
\usepackage{hyperref}

\newcommand{\dist}{\mathrm{dist}}

\newcommand{\set}[1]{\left\lbrace #1 \right\rbrace}
\newcommand{\norm}[1]{\left\Vert #1 \right\Vert}
\newcommand{\abs}[1]{\left\vert #1 \right\vert}

\theoremstyle{definition}
\newtheorem{theorem}{Theorem}[section]
\newtheorem{proposition}[theorem]{Proposition}
\newtheorem{corollary}[theorem]{Corollary}
\newtheorem{example}[theorem]{Example}
\newtheorem{definition}[theorem]{Definition}
\newtheorem{lemma}[theorem]{Lemma}
\newtheorem{definition-proposition}[theorem]{Definition-Proposition}
\newtheorem{remark}[theorem]{Remark}

\newtheorem{assumption}[theorem]{Assumption}

\begin{document}

\title{On funnel control with star-shaped funnels}

\author{Jonas Kirchhoff%
\thanks{J. Kirchhoff is with the Institute of Mathematics, Martin Luther University Halle-Wittenberg, Halle (Saale), Germany
        {\tt\small jonas.kirchhoff@mathematik.uni-halle.de}}%
}

\maketitle	
\thispagestyle{empty}
\pagestyle{empty}

\begin{abstract}
Funnel control for the trivial control-affine prototype is extended to the case of star-shaped funnels, both in the special case that the dynamical system is formulated in an Euclidean space and the general case of a Riemannian manifold.
\end{abstract}

\section{Introduction}

In this note, we consider reference tracking with transient guarantees for the trivial control-affine prototype
\begin{align}\label{eq:super_prototype}
\tfrac{\mathrm{d}}{\mathrm{d}t}x = f(x)+u,
\end{align}
where the state $x$ evolves on a manifold $M$ and $f\in\mathfrak{X}(M)$ is a global vector field. The reference trajectory is a known curve $x_{\mathrm{ref}}:[0,\infty)\to M$, and the transient guarantee, that we are aiming to give, take the form of some neighbourhood $\mathcal{F}(t)$ of $x_{\mathrm{ref}}(t)$ for each time $t\in[0,\infty)$. Precisely, the control aim is to obtain a feedback controller $u(t) = u(t,x(t))$ so that the closed loop of~\eqref{eq:super_prototype} admits, for each initial state $x(t_0)\in\mathcal{F}(t_0)$ a(n unique) solution defined on $[t_0,\infty)$ so that $x(t)\in\mathcal{F}(t)$ for all $t\in[t_0,\infty)$.

When $M = \mathbb{R}^n$ and each $\mathcal{F}(t)$ is an Euclidean ball $\mathcal{F}(t) = \mathbb{B}\big(x_{\mathrm{ref}}(t),\psi(t)\big)$ with sufficiently regular reference trajectory and funnel diameter $\psi:[0,\infty)\to\mathbb{R}_{>0}$, it is well-known that the control aim is satisfied by the \textit{funnel controller}, cf.~\cite{IlchRyanSang02,IlchRyan08,BergIlchRyan21,BergIlchRyan25}. In the companion paper~\cite{KircBerg26}, this was generalised to the case that $M$ is a smooth Riemannian manifold and each $\mathcal{F}(t)$ is a Riemannian ball with smooth funnel diameter. A generalisation of the classical funnel controller that removes the restriction to balls was proposed in~\cite{MehdLindBechDima23}, where funnels of the form
\begin{align*}
\mathcal{F}(t) = \bigcap_{i=1}^m\psi_i(t,\cdot)^{-1}(\set{(0,\infty)}),\qquad t\in[0,\infty)
\end{align*}
for some functions $\psi_1,\ldots,\psi_m:[0,\infty)\times\mathbb{R}^n\to\mathbb{R}$ were considered.

Here, we propose a further generalisation to funnel regions $\mathcal{F}(t)$, which are star-shaped with respect to the reference $x_{\mathrm{ref}}(t)$. The star-shape is the minimal requirement for a controller akin to the classical funnel controller, which considers $u$ as the scaled velocity along the geodesic from the current state to the current reference, to be applicable without any further information on the vector field and the reference trajectory, besides sufficient regularity.

The note is organised as follows: First, we consider the Euclidean case $M = \mathbb{R}^n$ and give a definition of star-shaped funnel regions. The well-posedness of the closed loop is proven in detail. In this proof, the sufficient conditions for global bounds on the distance to boundary and gain are identified. We then sketch out the generalisation to the case of a general Riemannian manifold and briefly discuss the advantages compared to Riemannian funnel control of~\cite{KircBerg26}.

\section{Euclidean funnel control with star-shaped funnel}

\subsection{Star-shaped moving domains}

\begin{definition}
A moving domain in $\mathbb{R}^n$ is a subset $\mathcal{F}(\cdot)\subseteq[0,\infty)\times\mathbb{R}^n$ so that, for all $t\in[0,\infty)$,
\begin{align*}
\mathcal{F}(t) := \set{x\in\mathbb{R}^n\,\big\vert\,(t,x)\in\mathcal{F}(\cdot)}
\end{align*}
is a non-empty open and connected set. $\mathcal{F}(\cdot)$ is star-shaped with respect to a curve $z:[0,\infty)\to\mathbb{R}^n$ if, and only if, each $\mathcal{F}(t)$ is star-shaped with respect to $z(t)$, i.e. $\set{z(t)+\lambda (x-z(t))\,\big\vert\,\lambda\in [0,1], x\in\mathcal{F}(t)} = \mathcal{F}(t)$.
\end{definition}

The prototypical example of a star-shaped moving domain are moving sublevel-sets of convex functions, and finite intersections of moving domains that are star-shaped with respect to the same curve.

\begin{example}
\begin{enumerate}[(i)]
\item Let $z:[0,\infty)\to\mathbb{R}^n$, and let $f:[0,\infty)\times\mathbb{R}^n\to\mathbb{R}$ so that each $f(t,\cdot)$, $t\in[0,\infty)$, is a continuous convex function. Choose further a curve $\psi:[0,\infty)\to\mathbb{R}$ so that $f(t,z(t))<\psi(t)$ for all $t\in[0,\infty)$. The set
\begin{small}
\begin{align*}
\mathcal{F}(\cdot;f,z,\psi) := \set{(t,x)\in[0,\infty)\times\mathbb{R}^n\,\big\vert\,f(t,x)<\psi(t)}
\end{align*}
\end{small}
is open, nonempty and convex, hence star-shaped with respect to $z$. The funnel domain of classical funnel control is given by $f(t,\cdot)\equiv\norm{\cdot-z(\cdot)}_2$, in which case we shall write $\mathcal{F}(\cdot;z,\psi)$
\item Let $z:[0,\infty)\to\mathbb{R}^n$ be a curve, and let $\mathcal{F}_1(\cdot),\ldots,\mathcal{F}_k(\cdot)$ be moving domains, each of whom is star-shaped with respect to $z$. Put $\mathcal{F}(\cdot) := \bigcap_{i=1}^k\mathcal{F}_i(\cdot)$ and observe that $z(t)\in \mathcal{F}(t) = \bigcap_{i=1}^k\mathcal{F}_i(t)$ is the non-empty finite intersection of open sets, and satisfies $z(t)+\lambda (x-z(t))\in\mathcal{F}(t)$ for all $\lambda\in[0,1]$ and $x\in\mathcal{F}(t)$. Therefore, $\mathcal{F}(\cdot)$ is a moving domain, which is star-shaped with respect to $z$. In this prototypical example, the intersection of this note and the work~\cite{MehdLindBechDima23} is to be found.
\end{enumerate}
\end{example}

\subsection{Regularity assumptions}

In classical Euclidean funnel control, the gain factor is the inverse distance to the funnel boundary, and the latter is the difference between funnel diameter and distance to the reference. In star-shaped Euclidean funnel control, not each time-slice of the funnel is assumed to be a ball in a (non-linear) metric and hence there is, in general, no funnel diameter. The distance to the boundary, however, does exist, and takes its place. A first observation is that the local Lipschitz continuity of the two functions is, in the classical case, equivalent.

\begin{proposition}
Let $z:[0,\infty)\to \mathbb{R}^n$ be locally Lipschitz continuous. A function $\psi:[0,\infty)\to(0,\infty)$ is locally Lipschitz continuous if, and only if, the function
\begin{align*}
\delta_{\mathcal{F}(\cdot;z,\psi)}:[0,\infty)\times\mathbb{R}^n & \to\mathbb{R},\\
(t,x) & \mapsto\mathrm{dist}\big(x,\mathbb{R}^n\setminus\mathcal{F}(t;\norm{\cdot},z,\psi)\big),
\end{align*}
where the distance is measured in the Euclidean norm, is locally Lipschitz continuous.
\end{proposition}
\begin{proof}
First, observe that local Lipschitz continuity only depends on the uniform equivalence class of the metric. Since all norms on $\mathbb{R}^{1+n}$ are equivalent, and the induced metrics of equivalent norms are uniform equivalent, we may without loss of generality equip $\mathbb{R}^{1+n}$ with the norm
\begin{align*}
\norm{(t,x)} := \abs{t}+\norm{x}_2,
\end{align*}
where $\norm{\cdot}_2$ denotes the Euclidean norm on $\mathbb{R}^n$; $[0,\infty)\times\mathbb{R}^n$ is equipped with the restriction of the induced metric. 

Furthermore, observe that, for all $t\in[0,\infty)$ and $x\in\mathbb{R}^n$,
\begin{align*}
\mathrm{dist}\big(x,\mathbb{R}^n\setminus\mathcal{F}(t;\norm{\cdot},z,\psi)\big) = \max\set{0,\psi(t)-\norm{x-z(t)}_2}.
\end{align*}
Let now $\delta_{\mathcal{F}(\cdot;z,\psi)}$ be locally Lipschitz continuous, and let $K\subseteq[0,\infty)$ be a compact set. Since $z$ is, in particular, continuous, there exists a compact set $K'\subseteq\mathbb{R}^n$ so that $z(t)\in K'$ for all $t\in K$. Let $L>0$ and $L'>0$ be the Lipschitz constants for $\delta_{\mathcal{F}(\cdot;z,\psi)}$ and $z$ on $K\times K'$ and $K$, respectively. Then, we have, for all $t,t'\in K$,
\begin{align*}
\abs{\psi(t)-\psi(t')} & = \abs{\delta_{\mathcal{F}(\cdot;z,\psi)}(t,z(t))-\delta_{\mathcal{F}(\cdot;z,\psi)}(t',z(t'))}\\
& \leq L\big(\abs{t-t'}+\norm{z(t)-z(t')}_2\big)\\
& \leq 2\max\set{L,LL'}\abs{t-t'}.
\end{align*}
Therefore, $\psi$ is indeed locally Lipschitz continuous.

Conversely, assume that $\psi$ is locally Lipschitz, and let $K\subseteq[0,\infty)\times \mathbb{R}^n$ be a compact set. Continuity of the Cartesian projection implies that
\begin{align*}
\widetilde{K}:=\set{t\in[0,\infty)\,\big\vert\,\exists x\in\mathbb{R}^n: (t,x)\in K}
\end{align*}
is compact. Therefore, $\psi$ is Lipschitz continuous on $\widetilde{K}$ with Lipschitz constant $L>0$. Let $L'$ be the Lipschitz constant for $z$ on $\widetilde{K}$. Let $(t,x),(t',x')\in K$, and assume first that $x\in\mathcal{F}(t;z,\psi)$ and $x'\in\mathcal{F}(t';z,\psi)$. Then,
\begin{align*}
& \abs{\delta_{\mathcal{F}(\cdot;z,\psi)}(t,x)-\delta_{\mathcal{F}(\cdot;z,\psi)}(t',x')}\\
& \qquad = \abs{\psi(t)-\psi(t')+\norm{x-z(t)}_2-\norm{x'-z(t')}_2}\\
& \qquad \leq \abs{\psi(t)-\psi(t')}+\abs{\norm{x-z(t)}-\norm{x'-z(t')}}\\
& \qquad \leq \abs{\psi(t)-\psi(t')}+\norm{x-x'+z(t)-z(t')}_2\\
& \qquad \leq \abs{\psi(t)-\psi(t')}+\norm{z(t)-z(t')}_2+\norm{x-x'}_2\\
& \qquad \leq \max\set{L,L',1}\big(\abs{t-t'}+\norm{x-x'}_2\big).
\end{align*}
In case that $x\in\mathcal{F}(t;\norm{\cdot},z,\psi)$ and $x'\notin\mathcal{F}(t';\norm{\cdot},z,\psi)$, we have
\begin{align*}
& \abs{\delta_{\mathcal{F}(\cdot;z,\psi)}(t,x)-\delta_{\mathcal{F}(\cdot;z,\psi)}(t',x')}\\
& \qquad = \psi(t)-\norm{x-z(t)}_2\\
& \qquad \leq \psi(t)-\norm{x-z(t)}_2-\big(\psi(t')-\norm{x'-z(t')}_2\big)\\
& \qquad = \abs{\psi(t)-\norm{x-z(t)}_2-\big(\psi(t')-\norm{x'-z(t')}_2\big)}\\
& \qquad \leq \max\set{L,L',1}\big(\abs{t-t'}+\norm{x-x'}_2\big);
\end{align*}
the case $x\notin\mathcal{F}(t;\norm{\cdot},z,\psi)$ and $x'\in\mathcal{F}(t';\norm{\cdot},z,\psi)$ is treated analogously. In the remaining case that both $x\notin\mathcal{F}(t;\norm{\cdot},z,\psi)$ and $x'\notin\mathcal{F}(t';\norm{\cdot},z,\psi)$, $\delta_{\mathcal{F}(\cdot;z,\psi)}(t,x) = \delta_{\mathcal{F}(\cdot;z,\psi)}(t',x') = 0$ and therefore we conclude that $\delta_{\mathcal{F}(\cdot;z,\psi)}$ is Lipschitz continuous on $K$ with Lipschitz constant $\max\set{L,L',1}$.
\end{proof}

Moreover, the classical funnel is automatically locally finite, which has to be assumed for star-shaped funnels; similar as it was assumed for Riemannian funnel control. This gives the following assumptions on the funnel domain and the reference trajectory.

\begin{assumption}\label{ass:Euclidean_assumption}
The moving domain $\mathcal{F}(\cdot)\subseteq[0,\infty)\times\mathbb{R}^n$ and the reference trajectory $x_{\mathrm{ref}}:[0,\infty)\to\mathbb{R}^n$ satisfy
\begin{enumerate}[(i)]
\item $\mathcal{F}(\cdot)$ is star-shaped with respect to $x_{\mathrm{ref}}$.
\item $\mathcal{F}(\cdot)$ is locally compactly bounded, i.e., for every $T\in(0,\infty)$, there exists a compact set $K(T)\subseteq[0,\infty)\times\mathbb{R}^n$ with $\mathcal{F}(\cdot)\vert_{[0,T]} := \mathcal{F}(\cdot)\cap [0,T]\times\mathbb{R}^n\subseteq K(T)$.
\item The function
\begin{equation}\label{eq:the_delta_function}
\begin{aligned}
\delta_{\mathcal{F}}:[0,\infty)\times\mathbb{R}^n & \to\mathbb{R},\\ \delta_{\mathcal{F}}(t,x) & := \dist\big(x,\mathbb{R}^n\setminus\mathcal{F}(t)\big),
\end{aligned}
\end{equation}
where the distance is measured in the Euclidean norm, is locally Lipschitz continuous.
\end{enumerate}
\end{assumption}

Since $\mathcal{F}(\cdot)$ is locally compactly bounded, each time-slice $\mathcal{F}(t)$, $t\in[0,\infty)$ is bounded. Therefore, $\mathbb{R}^n\setminus\mathcal{F}(t)$ is nonempty, and $\delta_{\mathcal{F}}$ defined in~\eqref{eq:the_delta_function} is well-defined. Moreover, $\delta_{\mathcal{F}}(t,\cdot)$ decreases along the rays towards the star-core of $\mathcal{F}(t)$ at most linearly.

\begin{lemma}\label{lem:ray-wise_superlipschitz}
Let $\mathcal{F}(\cdot)\subseteq[0,\infty)\times\mathbb{R}^n$ and $x_{\mathrm{ref}}:[0,\infty)\to\mathbb{R}^n$ satisfy Assumption~\ref{ass:Euclidean_assumption}. For each $t\in[0,\infty)$, $x\in\mathcal{F}(t)$ and $\lambda\in[0,1]$,
\begin{align}\label{eq:superlinearity}
\lambda\delta_\mathcal{F}(t,x) \leq \delta_{\mathcal{F}}\big(t,x_{\mathrm{ref}}(t)+\lambda(x-x_{\mathrm{ref}}(t))\big).
\end{align}
\end{lemma}
\begin{proof}
For $\lambda = 0$, nothing needs to be shown. Let $\lambda\in(0,1]$ and put $x_\lambda := x_{\mathrm{ref}}(t)+\lambda(x-x_{\mathrm{ref}}(t))\in\mathcal{F}(t)$ so that
\begin{align*}
x = \frac{1}{\lambda}\big(x_\lambda-(1-\lambda)x_{\mathrm{ref}}(t)\big).
\end{align*}
Let $y\in\mathbb{B}\big(x_\lambda,\lambda\delta_{\mathcal{F}}(t,x)\big)$ and define
\begin{align*}
y_\lambda := \frac{1}{\lambda}\big(y-(1-\lambda)x_{\mathrm{ref}}(t)\big).
\end{align*}
We observe that
\begin{align*}
\norm{y_\lambda-x}_2 = \frac{1}{\lambda}\norm{x_\lambda-y}_2 < \frac{1}{\lambda}\lambda\delta_{\mathcal{F}}(t,x)
\end{align*}
By definition of the distance function, we have, for all $z\in\mathcal{F}(t)$,
\begin{align}\label{eq:characterisation_delta}
\delta_{\mathcal{F}}(t,z) = \max\set{\delta>0\,\big\vert\,\mathbb{B}(z,\delta)\subseteq\mathcal{F}(t)},
\end{align}
and hence, in particular, $y_\lambda\in\mathcal{F}(t)$. Since $\mathcal{F}(t)$ is star-shaped with respect to $x_{\mathrm{ref}}(t)$
\begin{align*}
y = x_{\mathrm{ref}}(t)+\lambda (y_\lambda-x_{\mathrm{ref}}(t))\in\mathcal{F}(t)
\end{align*}
and hence $\mathbb{B}\big(x_\lambda,\lambda\delta_{\mathcal{F}}(t,x)\big)\subseteq\mathcal{F}(t)$. Applying the characterisation~\eqref{eq:characterisation_delta} once more gives that $\delta_{\mathcal{F}}(t,x_\lambda)\geq\lambda\delta_{\mathcal{F}}(t,x)$.
\end{proof}

Observe that the distance to the boundary may shrink along the ray from a point towards the reference trajectory, which is illustrated with an example.

\begin{example}
Consider the three points $A = (0,1)$, $B = (0,-1)$ and $C = (\sqrt{3},0)$, whose convex span is an equilateral triangle with interior
\begin{align*}
\Delta := \set{\alpha A+\beta B+\gamma C\,\big\vert\,\alpha,\beta,\gamma\in(0,1),\alpha+\beta+\gamma = 1}.
\end{align*}
In particular, $\Delta$ is convex and therefore star-shaped with respect to all of its points. Consider $x = \big(\frac{1}{\sqrt{3}},0\big)$, so that $\mathrm{dist}\big(x,\mathbb{R}^2\setminus\Delta\big) = \frac{1}{\sqrt{3}}$. With $z = (\varepsilon,0)$ for $0<\varepsilon<\frac{1}{\sqrt{3}}$, we observe
\begin{align*}
\mathrm{dist}\big(z+\lambda(x-z),\mathbb{R}^2\setminus\Delta\big) = \varepsilon+\lambda\left(\frac{1}{\sqrt{3}}-\varepsilon\right)
\end{align*}
for all $\lambda\in [0,1]$. In particular, the distance to boundary shrinks towards $z$.
\end{example}

Distance functions of co-bounded sets are Lipschitz continuous with Lipschitz constant 1. This gives the following.

\begin{lemma}\label{lem:slicewise_lipschitz}
Let $\mathcal{F}(\cdot)\subseteq[0,\infty)\times\mathbb{R}^n$ be a moving domain with $\mathcal{F}(t)\neq\mathbb{R}^n$, and consider $\delta_{\mathcal{F}}$ defined in~\eqref{eq:the_delta_function}. For each $t\in[0,\infty)$, $\delta_{\mathcal{F}}(t,\cdot)$ is Lipschitz-continuous with Lipschitz-constant 1.
\end{lemma}
\begin{proof}
Let $t\in[0,\infty)$, and let $x,y\in\mathcal{F}(t)$. Let $(\overline{x}_k)_{k\in\mathbb{N}}\in(\mathbb{R}^n\setminus\mathcal{F}(t))$ be sequence with $\delta_{\mathcal{F}}(t,x) = \lim_{k\to\infty}\norm{x-\overline{x}_k}$. By triangle inequality,
\begin{align*}
\delta_{\mathcal{F}}(t,y)\leq\norm{y-x}_2+\norm{x-\overline{x}_k}_2\overset{k\to\infty}{\to} \norm{y-x}_2+\delta_{\mathcal{F}}(t,x).
\end{align*}
Interchanging the r\^{o}les of $x$ and $y$ and combining the two inequalities implies
\begin{align}\label{eq:slicewise_lipschitz}
\abs{\delta_{\mathcal{F}}(t,y)-\delta_{\mathcal{F}}(t,x)}\leq \norm{x-y}.
\end{align}
If $x\notin\mathcal{F}(t)$ or $y\notin\mathcal{F}(t)$, then the inequality~\eqref{eq:slicewise_lipschitz} is trivially satisfied.
\end{proof}

\subsection{The star-shaped funnel controller}

We now give a proof of the trivial control-affine prototype of Euclidean funnel control with star-shaped funnels.

\begin{theorem}\label{thm:Euclidean_case}
Let $f\in\mathcal{C}^1(\mathbb{R}^n,\mathbb{R}^n)$, $x_{\mathrm{ref}}\in\mathcal{C}^1([0,\infty),\mathbb{R}^n)$, and let $\mathcal{F}(\cdot)$ be a moving domain so that Assumption~\ref{ass:Euclidean_assumption} is satisfied. The closed loop 
\begin{align*}
\tfrac{\mathrm{d}}{\mathrm{d}t}x = f(x)+\frac{1}{\delta_{\mathcal{F}}(t,x)}\big(x_{\mathrm{ref}}(t)-x\big)
\end{align*}
is forward complete on $\mathcal{F}(\cdot)$.
\end{theorem}
\begin{proof}
By definition of $\delta_{\mathcal{F}}$, $\delta_{\mathcal{F}}(t,x)>0$ for all $(t,x)\in\mathcal{F}(\cdot)$ and therefore the vector field of the closed loop,
\begin{align*}
F(t,x) := \begin{pmatrix}
1\\f(x)+\frac{1}{\delta_{\mathcal{F}}(t,x)}\big(x_{\mathrm{ref}}(t)-x\big)
\end{pmatrix} 
\end{align*}
is well-defined on $\mathcal{F}(\cdot)$. Since the inverse is smooth, and since the composition of Lipschitz continuous functions is Lipschitz, $\frac{1}{\delta_{\mathcal{F}}}$ is locally Lipschitz continuous. Therefore, $F$ is locally Lipschitz continuous and the classical theorem of Picard and Lindel\"of, cf.~\cite[Theorem VII.8.14]{AE_II}, implies that each initial datum $(t_0,x_0)\in\mathcal{F}(\cdot)$ admits an integral curve $(t,x):[t_0,\omega)\to\mathcal{F}(\cdot)$ where $\omega\in (0,\infty]$, $x(t_0) = x_0$ and $t$ is the identity, and all such integral curves coincide on their common domain. Let $(t,x)$ be maximal amongst all such integral curves. It remains to verify that $\omega = \infty$. We proceed in steps.

\textsc{Step 1:} We derive a dynamical inequality for the tracking error
\begin{align*}
e:[t_0,\omega)\to\mathbb{R},\qquad t\mapsto x(t)-x_{\mathrm{ref}}(t).
\end{align*}
Differentiating its squared Euclidean norm gives
\begin{align*}
& \tfrac{\mathrm{d}}{\mathrm{d}t}\frac{1}{2}\norm{e(t)}_2^2\\
&\quad = \left\langle e(t),f(x(t))-\frac{e(t)}{\delta_{\mathcal{F}}(t,x)}-\tfrac{\mathrm{d}}{\mathrm{d}t}x_{\mathrm{ref}}(t)\right\rangle_2\\
&\quad \leq \norm{e(t)}_2\left(\norm{f(x(t))}_2+\norm{\tfrac{\mathrm{d}}{\mathrm{d}t}x_{\mathrm{ref}}(t)}_2\right)-\frac{\norm{e(t)}_2^2}{\delta_{\mathcal{F}}(t,x)}.
\end{align*}
In particular, we conclude
\begin{align}\label{eq:nonzero_error_dynamical_bound}
\tfrac{\mathrm{d}}{\mathrm{d}t}\norm{e(t)}_2\leq \norm{f(x(t))}_2+\norm{\tfrac{\mathrm{d}}{\mathrm{d}t}x_{\mathrm{ref}}(t)}_2-\frac{\norm{e(t)}_2}{\delta_{\mathcal{F}}(t,x)}
\end{align}
for all $t\in(t_0,\omega)$ with $e(t)\neq 0$.

\textsc{Step 2:} Let $T\in(t_0,\omega)$, and let $K(T)\subseteq[0,\infty)\times\mathbb{R}^n$ be a compact set with $\mathcal{F}(\cdot)\vert_{[0,T]}\subseteq K(T)$, cf. Assumption~\ref{ass:Euclidean_assumption}\,(ii).  We calculate a lower bound for $\tfrac{\mathrm{d}}{\mathrm{d}t}\delta_{\mathcal{F}}(t,x(t))$ for almost all $t\in[t_0,T]$. Since $\delta_{\mathcal{F}}$ is Lipschitz-continuous, $\delta_{\mathcal{F}}\circ(t,x)$ is Lipschitz-continuous on the compact interval $[t_0,T]$ and hence absolutely continuous. Therefore, $\delta_{\mathcal{F}}\circ(t,x)$ is almost everywhere differentiable, cf.~\cite[Theorem 7.18]{Rudi87}, and the question is well-posed. Let $t\in(t_0,T)$ so that $\delta_{\mathcal{F}}\circ(t,x)$ is differentiable in $t$. By Assumption~\ref{ass:Euclidean_assumption}\,(iii), $\delta_{\mathcal{F}}$ is Lipschitz-continuous on $K(T)$ with Lipschitz constant $L(T)>0$, so that, for $h\in\big(0,\delta_{\mathcal{F}}(t,x(t))\big)$ with $\abs{h}<\min\set{T-t,t-t_0,\delta_{\mathcal{F}}(t,x(t))}$,
\begin{align*}
\delta_{\mathcal{F}}(t,x(t+h))-\delta_{\mathcal{F}}(t+h,x(t+h))\leq L(T)h.
\end{align*}
Define
\begin{align*}
\overline{x}(t,h) & := x_{\mathrm{ref}}(t)+\frac{\delta_{\mathcal{F}}(t,x(t))-h}{\delta_{\mathcal{F}}(t,x(t))}(x(t)-x_{\mathrm{ref}}(t))\\
& = x(t)+\frac{h}{\delta_{\mathcal{F}}(t,x(t))}\big(x_{\mathrm{ref}}(t)-x(t)\big)
\end{align*}
By Lemma~\ref{lem:slicewise_lipschitz}, we have
\begin{align*}
\abs{\delta_{\mathcal{F}}\left(t,\overline{x}(t,h)\right)-\delta_{\mathcal{F}}(t,x(t+h))} \leq \norm{\overline{x}(t,h)-x(t+h)},
\end{align*}
and Lemma~\ref{lem:ray-wise_superlipschitz} implies
\begin{align*}
\delta_{\mathcal{F}}\left(t,\overline{x}(t,h)\right) \geq \frac{\delta_{\mathcal{F}}(t,x(t))-h}{\delta_{\mathcal{F}}(t,x(t))}\delta_{\mathcal{F}}(t,x(t)) = \delta_{\mathcal{F}}(t,x(t))-h.
\end{align*}
Combining the two inequalities gives
\begin{align*}
\delta_{\mathcal{F}}(t,x(t))-\delta_{\mathcal{F}}(t,x(t+h))\leq - \norm{\overline{x}(t,h)-x(t+h)}-h
\end{align*}
and therefore
\begin{equation}\label{eq:1_gleichnug}
\begin{aligned}
& \delta_{\mathcal{F}}(t+h,x(t+h)) - \delta_{\mathcal{F}}(t,x(t))\\
&\qquad\geq -(1+L(T))h-\norm{\overline{x}(t,h)-x(t+h)}.
\end{aligned}
\end{equation}
First order Taylor approximation gives now
\begin{align*}
& x(t+h)-x(t)\\
&\quad = hf(x(t))+\frac{h}{\delta_{\mathcal{F}}(t,x(t))}\big(x_{\mathrm{ref}}(t)-x(t)\big)+o(h)
\end{align*}
and hence
\begin{align*}
\overline{x}(t,h)-x(t+h) = hf(x(t))+o(h).
\end{align*}
Plugging into~\eqref{eq:1_gleichnug} yields
\begin{align*}
& \delta_{\mathcal{F}}(t+h,x(t+h)) - \delta_{\mathcal{F}}(t,x(t))\\
& \qquad\geq-(1+L(T)+F(T))h+o(h),
\end{align*}
where
\begin{align*}
F(T) & := \max_{(t,x)\in K(T)}\norm{f(x)}_2.
\end{align*}
Dividing by $h\neq 0$ and taking the limit for $h\to 0$ gives
\begin{align}\label{eq:lower_bound_for_the_delta_derivative}
\tfrac{\mathrm{d}}{\mathrm{d}t} \delta_{\mathcal{F}}(t,x(t)) \geq -(1+L(T)+F(T)) =: -\check{\delta}(T).
\end{align}
\textsc{Step 3}: We calculate a lower bound for $\delta_{\mathcal{F}}\circ(t,x)$ on $[t_0,T]$. Consider
\begin{align*}
\delta_{\mathrm{ref}}(T) := \min_{t\in[t_0,T]}\delta_{\mathcal{F}}(t,x_{\mathrm{ref}}(t))>0
\end{align*}
and put
\begin{align*}
V(T) := \max_{t\in[t_0,T]}\norm{\tfrac{\mathrm{d}}{\mathrm{d}t}x_{\mathrm{ref}}(t)}.
\end{align*}
Assume that there is some $\tau\in[t_0,T]$ with
\begin{align}\label{eq:distance_is_dangerously_close}
\delta_{\mathcal{F}}(\tau,x(\tau))<\frac{\delta_{\mathrm{ref}}(T)}{2}
\end{align}
and define
\begin{align*}
t_1 := \max\set{t_0,\sup\set{t\in[t_0,\tau]\,\big\vert\,\delta_{\mathcal{F}}(t,x(t))\geq\frac{\delta_{\mathrm{ref}}(T)}{2}}}
\end{align*}
and
\begin{align*}
t_2 := \min\set{T,\inf\set{t\in[\tau,T]\,\big\vert\,\delta_{\mathcal{F}}(t,x(t))\geq\frac{\delta_{\mathrm{ref}}(T)}{2}}}.
\end{align*}
Continuity of $\delta_{\mathcal{F}}\circ (t,x)$ gives $t_1<\tau<t_2$, and
Lemma~\ref{lem:slicewise_lipschitz} implies, for all $t\in [t_1,t_2]$,
\begin{align}\label{eq:error_bound_from_boundary_bound}
\norm{e(t)}_2\geq \delta_{\mathcal{F}}(t,x_{\mathrm{ref}}(t))-\delta_{\mathcal{F}}(t,x(t))\geq \frac{\delta_{\mathrm{ref}}(T)}{2}.
\end{align}
In particular, $\norm{e}_2$ is differentiable on $(t_1,t_2)$, and the dynamical inequality~\eqref{eq:nonzero_error_dynamical_bound} implies
\begin{align*}
\tfrac{\mathrm{d}}{\mathrm{d}t}\norm{e(t)}_2 \leq V(T)+F(T)-\frac{\delta_{\mathrm{ref}}(T)}{2\delta(t,x(t))}.
\end{align*}
In particular, we may integrate and obtain
\begin{align*}
\int_{t_1}^{t_2}\frac{1}{\delta(t,x(t))}\,\mathrm{d}t & \leq \frac{2}{\delta_{\mathrm{ref}}(T)}\Big(\norm{e(t_1)}-\norm{e(t_2)}\\
& \qquad +(V(T)+F(T))(t_2-t_1)\Big)
\end{align*}
The latter can be bounded from above by
\begin{align*}
\Delta(T) & :=  \frac{2}{\delta_{\mathrm{ref}}(T)}\Big((V(T)+F(T))(T-t_0)\\
& \hspace{2cm} + \max_{(t,x),(t',x')\in K(T)} \norm{x-x'}\Big)
\end{align*}
Consider the function
\begin{equation}\label{eq:das_vieh}
\begin{aligned}
\varphi:[t_1,t_2] & \to\mathbb{R},\\
t & \mapsto \delta_{\mathcal{F}}(t,x(t))\exp{\check{\delta}(T)\int_{t_1}^t\frac{1}{\delta(\vartheta,x(\vartheta))}\,\mathrm{d}\vartheta},
\end{aligned}
\end{equation}
which is almost everywhere differentiable with derivative
\begin{small}
\begin{align*}
\tfrac{\mathrm{d}}{\mathrm{d}t}\varphi(t) = \Big(\tfrac{\mathrm{d}}{\mathrm{d}t}\delta_{\mathcal{F}}(t,x(t))+\check{\delta}(T)\Big)\exp{\check{\delta}(T)\int_{t_1}^t\frac{1}{\delta(\vartheta,x(\vartheta))}\,\mathrm{d}\vartheta}
\end{align*}
\end{small}
(almost everywhere), which is, by~\eqref{eq:lower_bound_for_the_delta_derivative}, (almost everywhere) non-negative. Therefore, $\varphi$ is non-decreasing and we conclude from $\delta_{\mathcal{F}}(t_1,x(t_1)) = \varphi(t_1)$
\begin{align*}
& \delta_{\mathcal{F}}(t,x(t))\\
& \qquad\geq \delta_{\mathcal{F}}(t_1,x(t_1))\exp\left(-{\check{\delta}(T)\int_{t_1}^{t_2}\frac{1}{\delta(\vartheta,x(\vartheta))}\,\mathrm{d}\vartheta}\right)\\
& \qquad\geq \delta_{\mathcal{F}}(t_1,x(t_1))\exp\left(-\check{\delta}(T)\Delta(T)\right)\\
& \qquad = \min\set{\frac{\delta_{\mathrm{ref}}(T)}{2},\delta_{\mathcal{F}}(t_0,x(t_0))}\exp\left(-\check{\delta}(T)\Delta(T)\right)
\end{align*}
for all $t\in[t_1,t_2]$. Observe that the lower bound is independent of the specifically chosen $t_1$ and $t_2$, so that we find
\begin{align*}
& \delta_{\mathcal{F}}(t,x(t))\\
& \qquad \geq \min\set{\frac{\delta_{\mathrm{ref}}(T)}{2},\delta_{\mathcal{F}}(t_0,x(t_0))}\exp\left(-\check{\delta}(T)\Delta(T)\right)
\end{align*}
for all $t\in[t_0,T]$.

\textsc{Step 4}: We conclude that $\omega = \infty$. Seeking a contradiction, assume that $\omega$ is finite. Let $K(\omega)\subseteq[0,\infty)\times\mathbb{R}^n$ be a compact set so that $\mathcal{F}(\cdot)\vert_{[0,\omega]}\subseteq K(\omega)$. Without loss of generality, we may assume that $K(T)\subseteq K(\omega)$ for all $T\in (t_0,\omega)$, which gives us the constants
\begin{align*}
F(\omega) & := \max_{(t,x)\in K(\omega)}\norm{f(x)}_2\geq F(T),\\\delta_{\mathrm{ref}}(\omega) & := \min_{t\in[t_0,\omega]}\delta_{\mathcal{F}}(t,x_{\mathrm{ref}}(t)) \leq \delta_{\mathrm{ref}}(T),\\
V(\omega) & := \max_{t\in[t_0,\omega]}\norm{\tfrac{\mathrm{d}}{\mathrm{d}t}x_{\mathrm{ref}}(t)}\geq V(T),\\
\Delta(\omega) & :=  \frac{2}{\delta_{\mathrm{ref}}(\omega)}\Big((V(\omega)+F(\omega))(\omega-t_0)\\
& \hspace{2cm} + \max_{(t,x),(t',x')\in K(\omega)} \norm{x-x'}\Big)\geq\Delta(T).
\end{align*}
Applying Step 3 to all $T\in[t_0,\omega)$ implies the uniform positive lower bound
\begin{align*}
& \delta_{\mathcal{F}}(t,x(t))\\
& \qquad \geq \min\set{\frac{\delta_{\mathrm{ref}}(\omega)}{2},\delta_{\mathcal{F}}(t_0,x(t_0))}\exp\left(-\check{\delta}(\omega)\Delta(\omega)\right)=:\varepsilon
\end{align*}
for all $t\in[t_0,\omega)$. Since $\delta_{\mathcal{F}}$ is, in particular, continuous, the set
\begin{align*}
& \set{(t,x)\in[t_0,\omega]\times\mathbb{R}^n\,\big\vert\,\delta_{\mathcal{F}}(t,x)\geq\varepsilon}\\
&\hspace*{3cm} = \delta_{\mathcal{F}}^{-1}\big([\varepsilon,\infty)\big)\cap\big([t_0,\omega]\times\mathbb{R}^n\big)
\end{align*}
is closed as the intersection of closed sets. By definition of $\delta_{\mathcal{F}}$, we observe $\mathcal{F}(\cdot) = \delta_{\mathcal{F}}^{-1}\big((0,\infty)\big)$ and therefore 
\begin{align*}
\set{(t,x)\in[t_0,\omega]\times\mathbb{R}^n\,\big\vert\,\delta_{\mathcal{F}}(t,x)\geq\varepsilon}\subseteq\mathcal{F}(\cdot)\vert_{[0,\omega]}\subseteq K(\omega)
\end{align*}
is a closed subset of a compact set, hence compact. Since $\omega$ is maximal, there must be some $\widehat{\vartheta}\in[t_0,\omega)$ so that
\begin{align*}
(\widehat{\vartheta},x(\widehat{\vartheta}))\notin\set{(t,x)\in[t_0,\omega]\times\mathbb{R}^n\,\big\vert\,\delta_{\mathcal{F}}(t,x)\geq\varepsilon}
\end{align*}
or, equivalently, $\delta_{\mathcal{F}}(\widehat{\vartheta},x(\widehat{\vartheta}))<\varepsilon$, which contradicts the construction of $\varepsilon$.
\end{proof}

As an immediate consequence, we obtain the sufficient conditions for global bounds for $\delta_{\mathcal{F}}\circ(t,x)$ and $u\circ(t,x)$.

\begin{corollary}\label{cor:global_bounds}
Let $f\in\mathcal{C}^1(\mathbb{R}^n,\mathbb{R}^n)$ and $x_{\mathrm{ref}}\in\mathcal{C}^1([0,\infty),\mathbb{R}^n)$, and let $\mathcal{F}(\cdot)$ be a moving domain so that Assumption~\ref{ass:Euclidean_assumption} is satisfied. Consider $x\in\mathcal{C}^1([t_0,\infty),\mathbb{R}^n)$, $t_0\geq 0$, so that
\begin{align*}
\tfrac{\mathrm{d}}{\mathrm{d}t}x = f(x)+\frac{1}{\delta_{\mathcal{F}}(t,x)}\big(x_{\mathrm{ref}}(t)-x\big).
\end{align*}
If 
\begin{enumerate}[(i)]
\item $\delta_{\mathcal{F}}$ is globally Lipschitz with Lipschitz constant $L$,
\item $F := \sup_{(t,z)\in\mathcal{F}(\cdot)}\norm{f(z)} < \infty$,
\item $V := \sup_{t\in[0,\infty)}\norm{\tfrac{\mathrm{d}}{\mathrm{d}t}x_{\mathrm{ref}}(t)} < \infty$,
\item $\widehat{\delta} := \sup_{(t,z)\in\mathcal{F}(\cdot)}\delta_{\mathcal{F}}(t,z) < \infty$,
\item $\delta_{\mathrm{ref}} := \inf_{t\in[0,\infty)}\delta_{\mathcal{F}}(t,x_{\mathrm{ref}}(t)) > 0$,
\end{enumerate}
then
$\delta_{\mathcal{F}}(\cdot,x(\cdot))$ is bounded from below and $x_{\mathrm{ref}}-x$ is bounded.
\end{corollary}
\begin{proof}
Let $e := x_{\mathrm{ref}}-x$. The dynamical inequality~\eqref{eq:nonzero_error_dynamical_bound} implies
\begin{align*}
\tfrac{\mathrm{d}}{\mathrm{d}t}\norm{e}_2\leq F+V-\frac{\norm{e}_2}{\widehat{\delta}}<0
\end{align*}
whenever $e>\widehat{\delta}(F+V)$ and hence we observe
\begin{align}\label{eq:global_error_bound}
\norm{e}_2\leq\max\set{\norm{e(0)}_2,\widehat{\delta}(F+V)} =: E
\end{align}
In step 3 in the proof of Theorem~\ref{thm:Euclidean_case}, replace $\frac{\delta_{\mathrm{ref}}(T)}{2}$ by $\widetilde{\delta} := \frac{\delta_{\mathrm{ref}}}{4(1+V+F)}$. Then, whenever $\delta_{\mathcal{F}}(t,e)<\widetilde{\delta}$,
\begin{align*}
\frac{\norm{e(t)}_2}{\delta(t,x(t))}\geq 2(1+V+F),
\end{align*}
so that the dynamical inequality~\eqref{eq:nonzero_error_dynamical_bound} implies 
\begin{align*}
\tfrac{\mathrm{d}}{\mathrm{d}t}\norm{e(t)}_2\leq -\frac{1}{2\delta_{\mathcal{F}}(t,x(t))}\norm{e(t)}_2
\end{align*}
 and integration gives
\begin{align*}
\int_{t_1}^{t_2} \frac{\norm{e(t)}_2}{\delta_{\mathcal{F}}(t,x(t))}\,\mathrm{d}t\leq 2(\norm{e(t_1)}_2-\norm{e(t_2)}_2)\leq 2E
\end{align*}
on each interval $[t_1,t_2]$ with $\delta_{\mathcal{F}}(t,e)<\widetilde{\delta}$ for all $t\in [t_1,t_2]$. Using the analogous function as to~\eqref{eq:das_vieh}, we may conclude that
\begin{align*}
\delta_{\mathcal{F}}(t,x(t))\geq\min\set{\delta_{\mathcal{F}}(t_0,x(t_0)),\widetilde{\delta}}\exp{-\frac{4(1+L+F)E}{\delta_{\mathrm{ref}}}}
\end{align*}
for all $t\in[t_0,\infty)$, which gives a global upper bound for the gain $\frac{1}{\delta_{\mathcal{F}}(t,x(t))}$.
\end{proof}

\begin{remark}
$\delta_{\mathcal{F}}(\cdot,x_{\mathrm{ref}}(\cdot))$ is bounded from below if, and only if, there exists some positive $\lambda>0$ so that $\mathbb{B}(x_{\mathrm{ref}}(t),\lambda)\subseteq\mathcal{F}(t)$ for all $t\in[0,\infty)$. While, Theorem~\ref{thm:Euclidean_case} allows for the funnel diameter of classical Euclidean funnel control to decline arbitrarily small, this comes, in general, with the cost of unbounded gain and input.
\end{remark}

We shall briefly comment on the difference of the proof of Theorem~\ref{thm:Euclidean_case} compared to the classical situation of funnel control, where the time-slices are Euclidean balls around the reference, i.e., $\mathcal{F}(t) = \mathbb{B}(x_{\mathrm{ref}}(t),\psi(t))$ for all $t\in[0,\infty)$ and the locally Lipschitz-continuous funnel diameter $\psi$.

\begin{remark}
For classical funnel control with funnel diameter $\psi$ and funnel domain $\mathcal{F}(\cdot;\norm{\cdot}_2,x_{\mathrm{ref}},\psi)$, we observed the distance to boundary within the funnel domain is given as
\begin{align*}
\delta_{\mathcal{F}}(t,x) = \psi(t)-\norm{x_{\mathrm{ref}}(t)-x}
\end{align*}
for all $t\in[0,\infty)$ and $x\in\mathcal{F}(t;\norm{\cdot}_2,x_{\mathrm{ref}},\psi)$. Therefore, $\delta_{\mathcal{F}}$ is bounded away from zero along a solution $x$ if, and only if, $\norm{x_{\mathrm{ref}}-x}$ is bounded away from $\psi$. This is proven by assuming that $\norm{x_{\mathrm{ref}}-x}$ gets sufficiently close to $\psi$ and using the dynamical inequality~\eqref{eq:nonzero_error_dynamical_bound}. For general star-shaped funnel domains, the distance to boundary does not necessarily admit a closed form, so that the dynamical inequality~\eqref{eq:nonzero_error_dynamical_bound} does not contain enough information. Nonetheless, the general idea of the proof remains the same: Consider the case that the distance to the boundary falls below a certain threshold and observe that there is a mechanism which provides a lower bound.
\end{remark}

\section{Riemannian star-shaped funnels}

In this section, we give a rough sketch on the extension of Theorem~\ref{thm:Euclidean_case} to the general case that $M$ is a Riemannian manifold with Riemannian metric $g$. Recall ever so briefly the definitions of Riemannian geometry; for a textbook account, we refer e.g.~to~\cite{Lang95,Lee18}.

\begin{definition}
Let $M$ be a smooth manifold. A Riemannian metric is a smooth symmetric tensor field $g\in\Gamma(\mathcal{T}^*M\vee\mathcal{T}^*M)$ which is fibrewise positive definite, i.e. $\norm{X}_g^2 := \langle g_x,X\otimes X\rangle>0$ for all $x\in M$ and $X\in\mathcal{T}_xM\setminus\set{0}$. The r\^{o}le of shortest paths with constant velocity are played by geodesics, which are formally defined by means of the induced Levi-Civita connection, and the geodesics give rise to the Riemannian exponential map
\begin{align*}
\exp_{x}:\mathrm{dom}\,\exp_x\subseteq\mathcal{T}_xM\to M,
\end{align*}
which maps a tangent vector $X\in\mathcal{T}_xM$ to the value of the geodesic starting at $x$ with constant velocity $X$ at time $1$, whenever defined. The metric $g$ induces the Riemannian distance $d_g$ defined by integrating $\norm{\tfrac{\mathrm{d}}{\mathrm{d}t}\gamma(\cdot)}$ for piecewise smooth paths connecting two points in finite time and taking the infimum over all such paths.
\end{definition}

The Euclidean space is a particular example of a Riemannian manifold whose geodesics are straight lines and whose exponential map is the globally defined  and reads (in coordinates) $\exp_{x}X = x+X$. The direction $x_{\mathrm{ref}}-x$ is therefore $\exp_x^{-1}(x_{\mathrm{ref}})$, and the star-shape should be with respect to the geodesics of the metric. We obtain the following generalisation of Assumption~\ref{ass:Euclidean_assumption}.

\begin{assumption}\label{ass:Riemannian}
The funnel domain $\mathcal{F}(\cdot)\subseteq[0,\infty)\times M$ and the reference trajectory $x_{\mathrm{ref}}:[0,\infty)\to M$ satisfy
\begin{enumerate}[(i)]
\item For each $t\in[0,\infty)$, $\mathcal{F}(t):=\mathcal{F}(\cdot)\cap(\set{t}\times M)$ is open and contains $x_{\mathrm{ref}}(t)$.
\item For each $t\in[0,\infty)$ and $x\in\mathcal{F}(t)$, there exists a unique geodesic $\gamma:[0,1]\to M$ with $\gamma(0) = x$, $\gamma(1) = x_{\mathrm{ref}}(t)$ and $\gamma(\tau)\in\mathcal{F}(\tau)$ for all $\tau\in[0,1]$.
\item $\mathcal{F}(\cdot)$ is locally compactly bounded.
\item The function
\begin{equation}\label{eq:the_garbage_delta_function}
\begin{aligned}
\delta_{\mathcal{F}}:[0,\infty)\times M & \to\mathbb{R},\\ \delta_{\mathcal{F}}(t,x) & := \dist\big(x,M\setminus\mathcal{F}(t)\big),
\end{aligned}
\end{equation}
where the distance is measured in the Riemannian distance, is locally Lipschitz continuous.
\end{enumerate}
\end{assumption}

We may now formulate the Riemannian version of Theorem~\ref{thm:Euclidean_case} and, due to space constraints, give an ever so short sketch of the proof.

\begin{theorem}
Let $f\in\mathcal{X}(M)$ be a global $\mathcal{C}^1$-vector field, $x_{\mathrm{ref}}\in\mathcal{C}^1([0,\infty),M)$, and let $\mathcal{F}(\cdot)\subseteq[0,\infty)\times M$ be a moving domain so that Assumption~\ref{ass:Riemannian} is satisfied. The closed loop 
\begin{align}\label{eq:Riemannian_closed_loop}
\tfrac{\mathrm{d}}{\mathrm{d}t}x = f(x)+\frac{1}{\delta_{\mathcal{F}}(t,x)}\exp_{x}^{-1}(x_{\mathrm{ref}}(t))
\end{align}
is forward complete on $\mathcal{F}(\cdot)$.
\end{theorem}
\begin{proof}
By assumption~\ref{ass:Riemannian}\,(ii), for each $t\in[0,\infty)$ and $x\in\mathcal{F}(t)$, there exists a unique geodesic $\gamma$ from $x$ to $x_{\mathrm{ref}}(t)$ so that $\tfrac{\mathrm{d}}{\mathrm{d}t}\gamma(0) = \exp_x^{-1}(x_\mathrm{ref}(t))$. Therefore, the closed loop~\eqref{eq:Riemannian_closed_loop} is well-defined. Moreover, assumption~\ref{ass:Riemannian}\,(ii) implies that $\exp_{x_{\mathrm{ref}}(t)}$ maps injectively into $\mathcal{F}(t)$ and hence~\cite[Theorems 3.4 and 5.1]{Warn65} and~\cite[Theorem 18.1]{Miln63} together with the implicit function theorem give that the map $(t,x)\mapsto\exp_{x}^{-1}(x_{\mathrm{ref}}(t))$ is continuously differentiable on $\mathcal{F}(\cdot)$. Hence,~\eqref{eq:Riemannian_closed_loop} is a locally Lipschitz-continuous vector field, and each initial datum $(t_0,x_0)$ admits a unique forward-maximal solution $x\in\mathcal{C}^1([t_0,\omega),M)$. To show that $\omega = \infty$, one must replace Euclidean distances between points by the Riemannian distance between points, i.e. $\norm{\exp_{x}^{-1}(x_{\mathrm{ref}}(t))}_g$ takes the r\^{o}le of $\norm{e(t)}_2$, and the Euclidean norm between vectors by the fibrewise Riemannian norm, i.e. $\norm{f(x)}_2$ becomes $\norm{f(x)}_g$ and $\norm{\tfrac{\mathrm{d}}{\mathrm{d}t}x_{\mathrm{ref}}(t)}_2$ is replaced by $\norm{\tfrac{\mathrm{d}}{\mathrm{d}t}x_{\mathrm{ref}}(t)}_g$. For the calculation of the derivative (almost everywhere) in Step 2, replace $x_{\mathrm{ref}}(t)+\lambda(x(t)-x_{\mathrm{ref}}(t))$ by the analogous point on the connecting geodesic, $\exp_{x}\big((1-\lambda)\exp_x^{-1}(x_{\mathrm{ref}}(t))$ and observe that $\delta_{\mathcal{F}}(t,\exp_{x}\big((1-\lambda)\exp_x^{-1}(x_{\mathrm{ref}}(t))\big)\geq\lambda^{H(T)}\delta_{\mathcal{F}}(t,x(t))$, where $H(T) = \sup_{(t,x)\in K(T)}\norm{\nabla_x^2\norm{\exp_\cdot(x_{\mathrm{ref}}(t))}_g}$, which exists by the Hessian comparison theorem~\cite[Theorem 11.7]{Lee18} and where the norm is fibrewise the norm of a quadratic form.
\end{proof}

\section{Conclusion and outlook}

We have proposed a version of funnel control with star-shaped funnel domains. In the case of the Euclidean space, we have given the detailed proof, while the case of Riemannnian manifolds was only very roughly sketched out. A fun consequence of the Euclidean case is that the distance from the reference in Euclidean funnel control may be measured in any norm. In~\cite{KircBerg26}, this was observed for Hilbert-space norms, only.

In a follow-up paper, the results shall be generalised to the non-trivial control-affine prototype with arbitrary input bundle, possibly higher relative degree and nonlinear output.

\section*{Acknowledgement}

J.K. thanks Prof. Thomas Berger (Halle) for the constructive discussions.

\bibliographystyle{plain}
\bibliography{../Literatur}

\end{document}